# Reliable design of an integrated supply chain with expedited shipments under disruption risks

Jianxun Cui[a], Meng Zhao[a], Xiaopeng Li[1b], Mohsen Parsafard[b], Shi An[a]

a. School of Transportation Science and Engineering, Harbin Institute of Technology, Harbin, China 150001

b. Department of Civil and Environmental Engineering, University of South Florida, FL, USA 33620

ABSTRACT

This paper proposes a mathematical model for the design of a two-echelon supply chain where a set of suppliers serve a set of terminal facilities that receive uncertain customer demands. This model integrates a number of system decisions in both the planning and operational stages, including facility location, multi-level service assignments, multi-modal transportation configuration, and inventory management. In particular, we consider probabilistic supplier disruptions that may halt product supply from certain suppliers from time to time. To mitigate the impact from supplier disruptions, a terminal facility may be assigned to multiple supplies that back up each other with different service priorities. With such multi-level service assignments, even when some suppliers are disrupted, a facility may still be supplied by the remaining functioning suppliers. Expensive expedited shipments yet with assured fast delivery may be used in complement to less expensive regular shipments yet with uncertain long lead times. Combining these two shipment modes can better leverage the inventory management against uncertain demands. We formulate this problem into a mix-integer nonlinear program that simultaneously determines all these interdependent system decisions to minimize the expected system cost under uncertainties from both suppliers and demands. A customized solution algorithm based on the Lagrangian relaxation is developed to efficiently solve this model. Several numerical examples are conduced to test the proposed model and draw managerial insights into how the key parameters affect the optimal system design.

## 1. Introduction

Supply chain operations are susceptible to various uncertainties such as facility disruptions, transportation delays, and customer demand fluctuations. As evidenced in recent catastrophic events (e.g., West Coast Lockdown (Gibson, Defee, and Ishfaq 2015), Szechuan Earthquake (Chan, 2008), Fukushima nuclear leak (Holt, Campbell, and Nikitin 2012), Hurricane Sandy (Blake et al. 2013)), supply chain facilities are vulnerable to various natural and anthropogenic disruption risks such as floods, earthquakes, power outages, and labor actions. Such disruptions, once happening, can choke the supply of corresponding commodities (or services) at the very source. Even if the commodities are successfully sent out from the supply facilities, they may experience extensive transportation delays, especially when they are shipped with slow transportation modes (e.g., waterways and railroads (Tseng et al. 2005; Ouyang and Li 2010)). Such transportation delay may cause depletion of downstream stocks and loss of customer demand, particularly when customer demand is stochastic and fluctuating. To ensure customer service levels, one way is to hold a high inventory of commodities at the downstream terminals (or retailer stores), which however incurs excessive inventory holding cost. Or expedited transportation can be used to largely reduce the delivery time to avoid accumulation of unmet demand, which however may dramatically increase transportation cost due to expensive expedited services. Li (2013) showed that a better way would be wisely combining inventory management and expedited transportation such that neither a high inventory nor frequent expedited services are needed. This series of uncertainties throughout these interdependent planning and operational stages, if not properly managed, may seriously damage system performance and deteriorate customer satisfaction. An integrated design methodology is needed to plan an efficient and reliable supply chain system that not only smartly balances cost tradeoffs over space and time but also robustly hedges against the unexpected uncertainties from supply, transportation, inventory, and demand.

---

[1] Corresponding author. Tel: 813-974-0778, E-mail: xiaopengli@usf.edu.

There have been a number of studies addressing different facets of supply chain design. On the facility location side, one recently intensively investigated topic is reliable facility location design. Studies on this topic basically aim to increase the expected performance of a supply chain system across various facility disruption scenarios by adding proper redundancy to the location design. On the operational side, freight lead time uncertainties and customer demand fluctuations have been well recognized as major challenges to inventory management and customer service quality. A recent study by Li (2013) proposed a strategy integrating occasional expedited shipments and proper inventory management that can successfully reduce the expected system cost under these transportation and demand uncertainties. However, for a realistic supply chain system that faces both facility disruptions and operational uncertainties simultaneously, it is imperative to have a system design method that is not only robust against facility disruption risks but immune to operational uncertainties.

This paper aims to bridge this research gap by proposing an integrated supply chain system design model that simultaneously determines facility location, multi-modal transportation configuration, and inventory management decisions all together under both facility disruption risks and operational uncertainties. This model considers a two-echelon supply chain system where a set of downstream terminal facilities order products from a subset of candidate upstream suppliers per arriving customer demands. Each supplier is however subject to unexpected disruptions from time to time, which force the terminals that used to be served by this supplier to divert to other suppliers or completely lose the service. To assure the service reliability, a terminal may be assigned to a sequence of suppliers such that if some of them are disrupted, the terminal can resort to the remaining according to the assignment priorities. Each shipment from a supplier to a terminal can be delivered via either a regular transportation mode that is cheap yet has a long and uncertain lead time or an expedited transportation mode that is much more expensive yet assures timely delivery. The adaption of expedited services also affects a terminal's inventory position and the corresponding inventory holding cost. The system design of this problem is very challenging. Not to mention the inherited NP-hardness of a location problem, the system has to face an extremely large number of possible supplier disruption scenarios exponential to the number of the suppliers. Further, the nested uncertainties from transportation delays and customer arrivals will complicate this problem even more. With our efforts, a compact polynomial-size mathematical programming model is proposed that integrates all these decisions components, including supplier location selections, supplier assignments to terminals, expedited transportation activation rules and inventory holding positions, so as to minimize the expected system cost from both location planning and operations under various uncertainties. The compact structure of this model formulation allows the development of an efficient Lagrangian relaxation algorithm that can efficiently solve this problem to a near-optimum solution. Numerical examples show that the proposed model can yield a supply chain system design that minimizes the impacts from probabilistic supplier disruptions and also leverages expedited shipments and inventory management to balance tradeoffs between transportation and inventory costs.

The rest of this paper is organized as follows. Section 2 reviews related literature. Section 3 formulates the proposed mathematical model for the integrated reliable design or the studied supply chain system. Section 4 develops a customized solution algorithm based on Lagrangian relaxation. Section 5 conducts numerical studies and discusses the experiment results. Section 7 concludes this paper and briefly discusses future research directions.

## 2. Literature review

Studies on facility location can be traced back to about a century ago (Weber 1929). Earlier location models focused on the single tradeoff between one-time facility investment and day-to-day transportation cost (see Daskin (1995) and Drezner (1995) for a review on these developments). These fundamental models were later extended in a number of directions that largely enriched the contents of facility location models. Spatially, the fundamental two-layer supply structures were extended to multi-layer (or multi-echelon) topologies (Şahin and Süral 2007). Temporally, single-period stationary operations were generalized to multi-period dynamic operations (Melo, Nickel, and Da Gama 2006). The system service was extended from a single commodity to multiple commodities that share the supply chain infrastructure (Klose and Drexl 2005). Direct transportation was extended to less-than-truck-load operations that involve vehicle routing decisions (Laporte 1987; Salhi and Petch 2007). Most of these models assume that all components of the supply chain system behave deterministically and their actions are fully predictable.

In reality, however, uncertainties exist almost ubiquitously throughout all components in a supply chain. Studies in 1980s (Daskin 1982; Daskin 1983; ReVelle and Hogan 1989; Batta, Dolan, and Krishnamurthy 1989) pointed out the need for facility redundancy under stochastic demand. Later studies (Lee, Padmanabhan, and Whang 1997; Ouyang and Daganzo 2006; Ouyang and Li 2010) further recognized that demand uncertainties cause serious challenges to inventory management when transportation takes long and uncertain lead times. To address this problem, facility location design has been integrated into inventory management to balance the tradeoff between spatial inventory distribution and transportation (Daskin, Coullard, and Shen 2002; Shen, Coullard, and Daskin 2003; Shu, Teo, and Shen 2005; Shen and Qi 2007; Snyder, Daskin, and Teo 2007; Qi, Shen, and Snyder 2010; Chen, Li, and Ouyang 2011). In addition, using faster transportation can alleviate the need for keeping high inventory, and thus expedited shipments can be adopted in the supply chain system to improve the overall system performance (Taghaboni-Dutta 2003; Huggins and Olsen 2003; Caggiano, Muckstadt, and Rappold 2006; Zhou and Chao 2010; Li 2013; Qi and Lee 2014). Li (2013) proposed a supply chain design framework that integrates location planning, inventory management, and expedited shipment configuration to mitigate the impacts from uncertain transportation and stochastic customer demands to the long-term supply chain performance.

Another major source of uncertainties in supply chain operations is unexpected facility disruptions, which was however largely overlooked in the facility location design literature in the last century. The catastrophic disasters in the recent years however resumed the recognition of the need for siting redundant facilities to hedge against disruption risks, and a number of modeling methods have been introduced for reliable location design under independent (Snyder and Daskin 2005; Cui, Ouyang, and Shen 2010; Qi, Shen, and Snyder 2010; Chen, Li, and Ouyang 2011; Li and Ouyang 2011; Li and Ouyang 2012; Li, Ouyang, and Peng 2013; Bai et al. 2015; Shishebori, Snyder, and Jabalameli 2014) and correlated (Li and Ouyang 2010; Berman and Krass 2011; Liberatore, Scaparra, and Daskin 2012; Lu, Ran, and Shen 2015) disruption risks.

This study aims to integrate the supply chain system design methods counteracting demand and transportation uncertainties with those addressing facility disruption risks. The model proposed by Li (2013) is extended to incorporate probabilistic facility disruptions. This extension is not trivial at all because the effect of facility disruptions and that from demand and transportation uncertainties are highly coupled. For example, disruptions of facilities will reduce candidate suppliers to customer terminals, which may in consequence increase transportation uncertainties and cumulate more unmet demand. Therefore, the extended problem is of much higher complexity, and substantial modeling efforts are needed to develop a comprehensive yet computationally-tractable model to solve this problem.

## 3. Model formulation

For the convenience of the readers, the mathematical notation is summarized in Table 1.

**Table 1** Notation List.

| Symbol | Description |
|---|---|
| $d_j$ | Demand rate at the terminal $j$ |
| $e_{ij}$ | Unit expedited shipment cost from supplier $i$ to terminal $j$ |
| $f_i$ | Fixed cost to install supplier $i$ |
| $h_j$ | Unit inventory holding cost at facility $j$ |
| $q$ | Supplier disruption probability for the regular service |
| $r_{ij}$ | Unit regular shipment cost from supplier $i$ to terminal $j$ |
| $t_{ij}$ | Expected regular shipment lead time from supplier $i$ to terminal $j$ |
| $L$ | Maximum assignment level |
| $P_{ij}(S_j)$ | Stock-out probability at terminal $j$ with base stock $S_j$ and regular supplier $i$ |
| $S_j$ | Base-stock position at terminal $j$ |
| $\bar{S}_j$ | Maximum allowable base-stock position at terminal $j$ |
| $X_i$ | Whether supplier $i$ is installed for service |
| $Y_{ijl}$ | Whether supplier $i$ provides regular service to terminal $j$ at assignment level l |
| $Z_{ij}$ | Whether supplier $i$ provides expedited service to terminal $j$ |
| **I** | Set of candidate suppliers, indexed by $i$ |
| **J** | Set of terminal facilities, indexed by $j$ |
| $\mathbf{L} = \{1, 2, \dots, L\}$ | Set of assignment levels, indexed by $l$ |
| $\mathbf{S}_j = \{1, 2, \dots, \bar{S}_j\}$ | Set of candidate base-stock positions |

Figure 1 illustrates the studied supply chain system, which includes set of terminal facilities denoted by **J** and a set of candidate suppliers denoted by **I**. Each terminal $j \in \mathbf{J}$ receives discrete demand for a certain commodity from a fixed pool of customers over time. We assume that at each terminal $j$, demand units arrive randomly with an expected rate of $d_j$. To feed the arriving demand, we assume that each terminal $j$ initially keeps a base-stock position $S_j \in \mathbf{S}_j := \{0,1,2,\dots,\bar{S}_j\}$ where $\bar{S}_j$ is a given capacity of the inventory at $j$, and the cost of holding one unit base stock per unit time is $h_j \geq 0$. This yields the system inventory cost as follows

$$C^{\mathrm{H}} := \sum_{j \in \mathbf{J}} h_j S_j \tag{1}$$

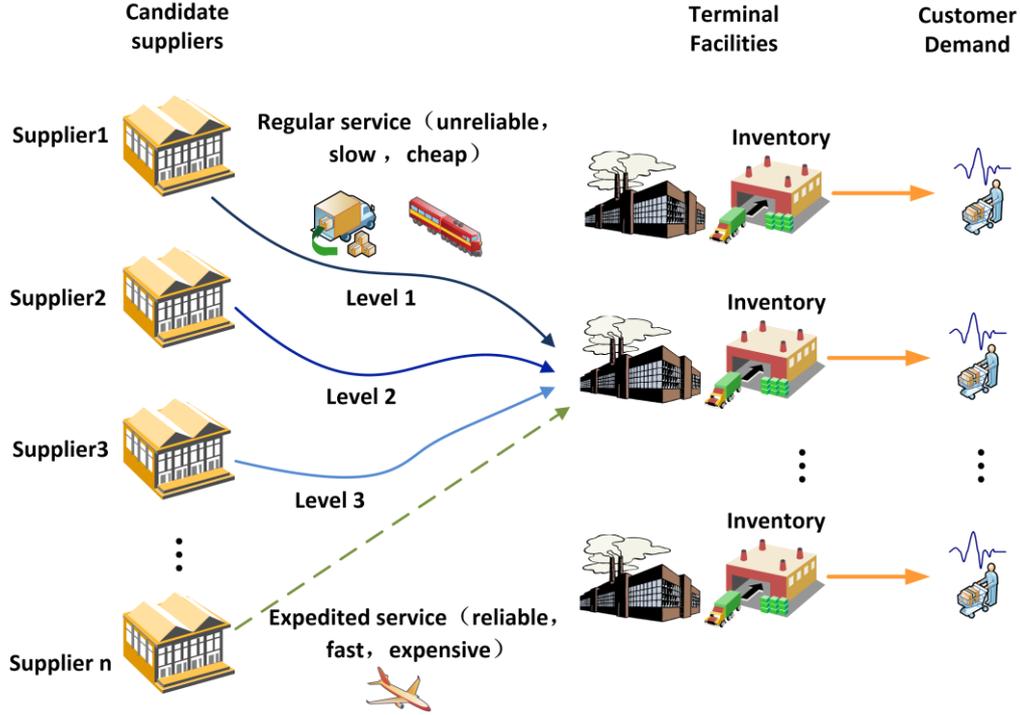

**Figure 1** Illustration of the studied supply chain system.

Whenever receiving a demand unit, terminal $j$ will first check its on-hand inventory and take one unit from this inventory, if any, to feed this demand unit. Meanwhile in order to maintain the base-stock inventory position, terminal $j$ can place an order right away from a supplier from **I**. This study considers possible supplier disruptions and assumes that each supplier can be disrupted independently any time at an identical probability $q$. To mitigate the impact from uncertain disruptions, a terminal is assigned to $L > 1$ suppliers at different priority levels for regular shipments. Every time, this terminal scans through these assigned suppliers from level 1 through level $L$ and places the order to the first functioning supplier. For the notation convenience, we define level set $\mathbf{L} := \{1,2,\dots,L\}$. In this way, the probability for a terminal to be served by its level-$l$ supplier is $(1-q)q^{l-1}, \forall l \in \mathbf{L}$. The assignments are specified by binary variables $[Y_{ijl}]_{i \in I, j \in J, l \in \mathbf{L}}$ such that $Y_{ijl} = 1$ if supplier $i$ is assigned to terminal $j$ at level $l$ or $Y_{ijl} = 0$ otherwise. Let $r_{ij}$ denote the cost to ship a unit commodity from supplier $i$ to terminal $j$, and then the total expected regular shipment cost is

$$C^{\mathrm{R}} := \sum_{i \in \mathbf{I}} \sum_{j \in \mathbf{J}} \sum_{l \in \mathbf{L}} d_j r_{ij} (1-q) q^{l-1} Y_{ijl}. \qquad (2)$$

We assume that the studied supply chain system has to maintain very high service quality such that customer demand has to be served right after it arrives. Despite being an economic option, a regular shipment is usually slow and unreliable. We assume that a regular shipment from supplier $i$ to terminal $j$ takes a random lead time with an expected value of $t_{ij}$. Since a longer shipment time is usually associated with a higher shipment cost for the same mode of transportation, we assume that $r_{ij} \geq r_{i'j} \Leftrightarrow t_{ij} \geq t_{i'j}, \forall i \neq i' \in \mathbf{I}, j \in \mathbf{J}$. In case that the regular shipments cannot arrive in time to meet the outstanding demand, the on-hand inventory at terminal $j$ may be depleted, particularly when the realized demand rate is high. In this case, terminal $j$ has to activate expedited transportation that always delivers shipments in a negligible lead time. We assume that every supplier provides an emergent expedited service that is independent from the regular service and never disrupts. When the on-hand inventory is depleted, a terminal uses the expedited service from a selected supplier, which however costs much more than regular transportation. Let $e_{ij}$ denote the cost to obtain an expedited shipment from supplier $i$ to terminal $j$, which shall satisfy $e_{ij} \gg r_{i'j}, \forall i' \in \mathbf{I}$. In order to quantify the expected expedited transportation cost, we will first quantify the probability for a terminal to activate the expedited service. Conditioning on that supplier $i$ is the active regular service provider to terminal $j$, the probability for terminal $j$ to use the expedited service can be represented as a function of initial inventory $S_j$ based on a truncated Poisson distribution (Li 2013),

$$P_{ij}(S_j) = \frac{(d_j t_{ij})^{S_j} / S_j!}{\sum_{s=0}^{S_j} (d_j t_{ij})^s / s!}. \qquad (3)$$

Note that once terminal $j$ places an expedited order from supplier $i'$, then no regular order is placed to the incumbent regular supplier $i$, and thus the actual additional cost due to this expedited order is $e_{i'j} - r_{ij}$. Define variables $[Z_{ij}]_{i \in I, j \in J}$ to denote the expedited service assignments such that $Z_{ij} = 1$ if terminal $j$'s expedited service provider is supplier $i$ or $Z_{ij} = 0$ otherwise. Then the total expected additional cost due to expedited shipments (or the marginal expedited cost) can be formulated as

$$C^{\mathrm{M}} := \sum_{j \in \mathbf{J}} \sum_{i \in \mathbf{I}} \sum_{l \in \mathbf{L}} \sum_{i' \in \mathbf{I}} d_j (e_{i'j} - r_{ij}) (1-q) q^{l-1} P_{ij}(S_j) Z_{i'j} Y_{ijl}. \qquad (4)$$

Another risk that the regular service is subject to is that all its suppliers may be disrupted simultaneously at probability $q^L$. If this happens, the regular service to this terminal becomes inactive, and we assume that it is now only served by emergency shipments from the previously assigned expedited supplier. The emergency cost structure stays the same as the previously defined expedited cost structure since they come from the same sources. Thus the expected system emergency cost is formulated

as

$$C^{\mathrm{E}} := \sum_{j \in \mathbf{J}} \sum_{i \in \mathbf{I}} d_j e_{ij} q^L Z_{ij}. \tag{5}$$

Finally, in this supply chain system, if candidate supplier $i$ is used by one or more terminals for either regular or expedited service, a fixed installation cost $f_i$ (prorated per unit time) is incurred. Define binary variables $[X_i]_{i \in \mathbf{I}}$ to denote the supplier location decisions such that $X_i = 1$ if candidate supplier $i$ is installed or $X_i = 0$ otherwise. This results in the system fixed installation cost as follows,

$$C^{\mathrm{F}} := \sum_{i \in \mathbf{I}} f_i X_i. \tag{6}$$

The system design includes integrated decisions of supplier location $[X_i]$, regular service assignments $[Y_{ijl}]$, expedited service assignments $[Z_{ij}]$, and initial inventory positions $[S_j]$ that collectively minimize the total system cost composed of (1), (2), (4), (5) and (6). Note that these cost components shall generally exhibit the following tradeoffs. Increasing supplier installations shall raise one-time fixed cost (6) but reduce day-to-day operational costs (2), (4) and (5). The higher inventory positions $[S_j]$ we set, which though increase inventory cost (1), the less frequent expedited shipments are needed according to probability function (3), and thus the less extra expedited transportation cost (4) is consumed. In order to quantitatively solve the detailed system design, the follow integer programming model is formulated.

$$\min C := \sum_{i \in \mathbf{I}} f_i X_i + \sum_{j \in \mathbf{J}} h_j S_j$$
$$+ \sum_{j \in \mathbf{J}} d_j \left\{ \sum_{i \in \mathbf{I}} \sum_{l \in \mathbf{L}} \left[ r_{ij} + \sum_{i' \in \mathbf{I}} (e_{i'j} - r_{ij}) Z_{i'j} P_{ij}(S_j) \right] (1-q) q^{l-1} Y_{ijl} + q^L \sum_{i \in \mathbf{I}} e_{ij} Z_{ij} \right\}, \tag{7}$$

$$\text{s.t.} \quad \sum_{l \in \mathbf{L}} Y_{ijl} - X_i \le 0, \quad \forall i \in \mathbf{I}, \ j \in \mathbf{J}, \tag{8}$$

$$Z_{ij} - X_i \le 0, \quad \forall i \in \mathbf{I}, \ j \in \mathbf{J}, \tag{9}$$

$$\sum_{i \in \mathbf{I}} Y_{ijl} = 1, \quad \forall j \in \mathbf{J}, \ l \in \mathbf{L}, \tag{10}$$

$$\sum_{i \in \mathbf{I}} Z_{ij} = 1, \quad \forall j \in \mathbf{J}, \tag{11}$$

$$S_j \in \mathbf{S}_j, \quad \forall j \in \mathbf{J}, \tag{12}$$

$$Y_{ijl} \in \{0,1\}, \quad \forall i \in \mathbf{I}, \ j \in \mathbf{J}, \ l \in \mathbf{L}, \tag{13}$$

$$Z_{ij} \in \{0,1\}, \quad \forall i \in \mathbf{I}, \ j \in \mathbf{J}, \tag{14}$$

$$X_i \in \{0,1\}, \quad \forall i \in \mathbf{I}. \tag{15}$$

Objective aims to minimize the summation of all cost components (2), (4), (5) and (6) across the entire system. Constraints (8) and (9) indicate that a supplier need to be installed first prior to its usage. Constraint (8) also ensures that if one of the suppliers is selected to provide the regular service to a terminal, it can only serve this terminal at one assignment level. Constraint (10) requires that one terminal has one and only one regular supplier at each level. Constraint (11) postulates that each terminal is assigned to one and only one expedited supplier. Constraints (12) - (15) are the corresponding integer and variable constraints for all variables.

**4. Solution approach**

Note that problem (7) - (15) is a complex nonlinear integer problem, which is apparently NP-hard because the basic incapacitated facility location problem is its special case. It will be extremely difficult to solve the exact optimum to a moderate-sized instance of this problem with off-the-shelf solvers. To tackle this challenge, this section proposes a customized solution approach that can solve a near-optimum solution to this problem very efficiently. Section 4.1 proposes a Lagrangian relaxation (LR) algorithm to obtain a lower bound to the optimal value of objective (7). Basically, the LR algorithm decomposes the original problem into a set of sub-problems that can be easily solved by simple enumeration. This relax solution however is likely infeasible. Section 4.2 proposes a heuristic algorithm to transform the relaxed solution to a feasible solution, which also provides an upper bound to the optimal objective(7). Finally, Section 4.3 adopts a sub-gradient search algorithm to iteratively update both upper and lower bound solutions to reduce the optimality gap between them until within an acceptable tolerance. Note that if the gap is reduced to zero, the final feasible solution is the true optimum. Otherwise, this feasible solution is a near-optimum solution no greater than the true optimum by the optimality gap.

**4.1 Lagrangian relaxation**

The Lagrangian relaxation algorithm basically relaxes constraints (8) and (9), the add the product of the-left hand side of these constraints and Lagrangian multipliers $\boldsymbol{\lambda} := \{\lambda_{ij} \geq 0\}_{i \in \mathbf{I}, j \in \mathbf{J}}$ and $\boldsymbol{\mu} := \{\mu_{ij} \geq 0\}_{i \in \mathbf{I}, j \in \mathbf{J}}$ into objective (7). We further add the following constraints

$$\sum_{l \in \mathbf{L}} Y_{ijl} \leq 1, \quad \forall i \in \mathbf{I}, \ j \in \mathbf{J}, \tag{16}$$

which are redundant to constraints (8) and are only used to improve the relaxed problem solution. This yields the following relax problem

$$\Delta(\boldsymbol{\lambda}, \boldsymbol{\mu}) := \min_{\mathbf{X},\mathbf{Y},\mathbf{Z},\mathbf{S}} \sum_{i \in \mathbf{I}} \left[ f_i - \sum_{j \in \mathbf{J}} \left( \lambda_{ij} + \mu_{ij} \right) \right] X_i \\ + \sum_{j \in \mathbf{J}} \left\{ \sum_{i \in \mathbf{I}} \sum_{l \in \mathbf{L}} \left( \sum_{i' \in \mathbf{I}} \alpha_{ii'jl} Z_{i'j} + \beta_{ijl} \right) Y_{ijl} + h_j S_j \right\}, \tag{17}$$

subject to (9)-(16), where

$$\alpha_{ii'jl} = d_j \left[ \left( e_{i'j} - r_{ij} \right) P_{ij} \left( S_j \right) (1-q) q^{l-1} + \frac{e_{i'j} q^L}{L} \right] + \frac{\mu_{i'j}}{L}, \tag{18}$$

and,

$$\beta_{ijl} = d_j r_{ij} (1-q) q^{l-1} + \lambda_{ij}. \tag{19}$$

Note in the above relaxed problem, the variables **X** are separated from the other variables. This allows us to decompose the relaxed problem into two sets of sub-problems. The first set only includes one sub-problem involving variables **X**:

$$\Gamma(\boldsymbol{\lambda}, \boldsymbol{\mu}) = \min_{\mathbf{X}} \sum_{i \in \mathbf{I}} \left[ f_i - \sum_{j \in \mathbf{J}} \left( \lambda_{ij} + \mu_{ij} \right) \right] X_i, \tag{20}$$

subject to binary constraint (15). Sub-problem (20) could be simply solved by setting $X_i = 1$ if $f_i - \sum_{j \in \mathbf{J}} (\lambda_{ij} + \mu_{ij}) \leq 0$ or $X_i = 0$ otherwise, which only takes a time complexity of $O(|\mathbf{I}||\mathbf{J}|)$. The second set contains $|\mathbf{J}|$ sub-problems, each associated with a terminal $j \in \mathbf{J}$, as follows:

$$\Phi_j(\boldsymbol{\lambda}, \boldsymbol{\mu}) := \min_{\{Y_{ijl}, Z_{ij}\}_{i \in \mathbf{I}, l \in \mathbf{L}}, S_j} \sum_{i \in \mathbf{I}} \sum_{l \in \mathbf{L}} \left( \sum_{i' \in \mathbf{I}} \alpha_{ii'jl} Z_{i'j} + \beta_{ijl} \right) Y_{ijl} + h_j S_j, \forall j \in \mathbf{J}, \tag{21}$$

subject to (10)-(14), (16) (where $\alpha_{ii'jl}$ and $\beta_{ijl}$ were formulated in (18) and (19)). We reformulate sub-problem (21) as a combinatorial problem to facilitate the solution algorithm. Define set $\mathbf{K} = \{(i_1, i_2, \cdots, i_L) | i_1 \neq i_2 \neq \cdots \neq i_L \in \mathbf{I}\}$, where each $(i_1, i_2, \cdots, i_L)$ specifies a strategy to assign the regular suppliers to terminal $j$ at all $L$ levels sequentially; i.e., supplier $i_l$ is assigned to terminal $j$ at level $l, \forall l \in \mathbf{L}$. For short we denote vector $(i_1, i_2, \cdots, i_L)$ with alias $k$. Then sub-problems (21) can be rewritten as:

$$\Phi_j(\boldsymbol{\lambda}, \boldsymbol{\mu}) := \min_{i' \in \mathbf{I}, k \in \mathbf{K}, S_j \in \mathbf{S}_j} C_{ki'j}(S_j) := A_{ki'j}(S_j) + h_j S_j + B_{ki'j}, \forall j \in \mathbf{J}, \tag{22}$$

where

$$A_{ki'j} = d_j \sum_{l \in \mathbf{L}} \left[ \left( e_{i'j} - r_{i_l j} \right) (1-q) q^{l-1} P_{i_l j}(S_j) \right], \tag{23}$$

$$B_{ki'j} = \sum_{l \in \mathbf{L}} \left[ d_j r_{ij} (1-q) q^{l-1} + \lambda_{ij} \right] + \frac{d_j e_{i'j} q^L + \mu_{i'j}}{L}. \tag{24} \quad \text{For given } k \text{ and } i',$$

$\min_{S_j \in \mathbf{S}_j} C_{ki'j}(S_j)$ can be solved with a bisection search method (BS) described in Appendix A. With this, problem (17) can be solved by a customized enumeration algorithm (EA) the does an exhaustive search through $k \in \mathbf{K}, i' \in \mathbf{I}$ for every $j \in \mathbf{J}$, as follows:

**Step EA1:** For each terminal $j \in \mathbf{J}$, we iterate through $(k, i') \in (\mathbf{K}, \mathbf{I})$ that specifies terminal $j$'s assignment strategy of both regular and expedited suppliers, and call the BS algorithm to solve $S_j^* := \arg\min_{S_j \in S} C_{ki'j}(S_j)$.

**Step EA2:** Find $(k^* = (i_1^*, i_2^*, \cdots, i_L^*), i'^*) := \arg\min_{k \in \mathbf{K}, i' \in \mathbf{I}} C_{ki'j}(S_j^*)$;

**Step EA3:** Return the optimal assignment strategy $(k^*, i'^*)$ and inventory position $S_j^*$;

**Step EA4:** Repeat EA1-3 for every supplier $j \in \mathbf{J}$ to get the optimal solution to $\mathbf{X}, \mathbf{Y}$ and $\mathbf{Z}$ as follows:

$$Y_{ijl} = \begin{cases} 1 & \text{if } i = i_l^*; \\ 0 & \text{otherwise,} \end{cases} \quad Z_{ij} = \begin{cases} 1 & \text{if } i = i'^*; \\ 0 & \text{otherwise,} \end{cases} \quad X_i = \max_{j \in \mathbf{J}, l \in \mathbf{L}} \{Y_{ijl}, Z_{ij}\}, \forall j \in \mathbf{J}, i \in \mathbf{I}, l \in \mathbf{L}. \quad (24)$$

Note that in the worst case, the time complexity of the EA algorithm for solving sub-problems (21) is $O(|\mathbf{J}||\mathbf{I}|^{L+1} \ln(\bar{S}_j))$.

By solving sub-problems (20) and (21), the object value of relaxed problem (17) for one set of given $\boldsymbol{\lambda}$ and $\boldsymbol{\mu}$ is equal to:

$$\Delta(\boldsymbol{\lambda}, \boldsymbol{\mu}) = \Gamma(\boldsymbol{\lambda}, \boldsymbol{\mu}) + \sum_{j \in \mathbf{J}} \Phi_j(\boldsymbol{\lambda}, \boldsymbol{\mu}), \quad (25)$$

which is a lower bound for the optimal value of problem (7)-(15) based on the duality property of Lagrangian relaxation (Geoffrion 1974). Note the time complexity of sub-problem (20) is of a lower order. Therefore, sub-problem (21) dominates the total time complexity of the relaxed problem (25).

**4.2 Feasible solutions**

If relaxed solution (24) happens to be feasible to the primal problem (7)-(15) with an identical objective value, then the solution will be optimal to the primal as well. Otherwise, if relaxed solution (24) is not feasible, which is likely for most large-scale instances, a feasible solution needs to be constructed by certain heuristic methods. One simple heuristic is to keep $\mathbf{Y}, \mathbf{Z}, \mathbf{S}$ values and adjust the $\mathbf{X}$ values as:

$$X_i = \max_{j \in \mathbf{J}, l \in \mathbf{L}} \{Y_{ijl}, Z_{ij}\}, \quad \forall i \in \mathbf{I}, \quad (26)$$

which could be solved very efficiently, i.e., in a time on the order of $O(|\mathbf{I}||\mathbf{J}|L)$. However, as $\mathbf{Y}, \mathbf{Z}$ values in the relaxed solution are usually much scattered, this feasible solution likely yields an excessive number of suppliers, leading to an unnecessarily high total cost. A better feasible solution is to fix $\mathbf{X}$ and adjust the other variables accordingly. Define $\bar{\mathbf{I}} := \{i | X_i = 1, \forall i \in \mathbf{I}\}$ to be the set of installed suppliers in the relaxed solution. Then, by setting $\lambda_{ij} = \mu_{ij} = 0$ and replacing $\mathbf{I}$ with $\bar{\mathbf{I}}$ in sub-problems(21), other feasible variable values can be determined by solving sub-problems (26)-(29):

$$\bar{\Phi}_j(\boldsymbol{\lambda}, \boldsymbol{\mu}) := \min_{\{Y_{ijl}, Z_{ij}\}_{i \in \bar{\mathbf{I}}, l \in \mathbf{L}}, S_j} \sum_{i \in \bar{\mathbf{I}}} \sum_{l \in \mathbf{L}} \left( \sum_{i' \in \bar{\mathbf{I}}} \bar{\alpha}_{ii'jl} Z_{i'j} + \bar{\beta}_{ijl} \right) Y_{ijl} + h_j S_j, \quad \forall j \in \mathbf{J}, \quad (26)$$

Subject to $\sum_{l \in \mathbf{L}} Y_{ijl} \leq 1, \quad \forall i \in \bar{\mathbf{I}}, \ j \in \mathbf{J},$ (27)

$$\sum_{i \in \bar{\mathbf{I}}} Y_{ijl} = 1, \quad \forall j \in \mathbf{J}, \ l \in \mathbf{L},$$ (28)

$$\sum_{i \in \bar{\mathbf{I}}} Z_{ij} = 1, \quad \forall j \in \mathbf{J},$$ (29)

where

$$\bar{\alpha}_{ii'jl} = d_j \left[ (e_{i'j} - r_{ij}) P_{ij}(S_j)(1-q)q^{l-1} + \frac{e_{i'j} q^L}{L} \right],$$ (30)

$$\bar{\beta}_{ijl} = d_j r_{ij} (1-q) q^{l-1}.$$ (31)

Then similar to the transformation from (21) to (22), we also define set $\bar{\mathbf{K}} = \{(i_1, i_2, \cdots, i_L) | i_1 \neq i_2 \neq \cdots \neq i_L \in \bar{\mathbf{I}}\}$ as all the strategies to assign regular suppliers from $\bar{\mathbf{I}}$ to terminal $j$ at all $L$ levels, and we also use alias $k$ to represent vector $(i_1, i_2, \cdots, i_L)$ for short. Then the transformed sub-problems are formulated as

$$\bar{\Phi}_j(\lambda, \mu) := \min_{i' \in \bar{\mathbf{I}}, k \in \bar{\mathbf{K}}, S_j \in \mathbf{S}_j} \bar{C}_{ki'j}(S_j) := \bar{A}_{kj} + \bar{B}_{ki'j} + h_j S_j, \forall j \in \mathbf{J}$$ (32)

where

$$\bar{A}_{kj} = d_j \sum_{l \in \mathbf{L}} r_{i_l j} (1-q) q^{l-1} \left(1 - P_{i_l j}(S_j)\right),$$ (33)

$$\bar{B}_{ki'j} = e_{i'j} \left( \frac{d_j q^L}{L} + d_j \sum_{l \in \mathbf{L}} (1-q) q^{l-1} P_{i_l j}(S_j) \right).$$ (34)

The exact optimal solution to each sub-problem (26)-(29) with any $j \in \mathbf{J}$ can be solved as follows. First denote $\bar{i}_j'^* := \operatorname{argmin}_{i' \in \bar{\mathbf{I}}} e_{i'j}$. Then we denote with vector $\bar{k}_j^* = \left( \bar{i}_1^{j*}, \bar{i}_2^{j*}, \cdots, \bar{i}_L^{j*} \right) \in \bar{\mathbf{K}}$ the first $L$ regular service suppliers sorted by the shipment cost to terminal $j$, i.e., $r_{\bar{i}_l^{j*} j} \leq r_{\bar{i}_m^{j*} j} \leq r_{ij}, \forall l < m \in \mathbf{L}, i \notin \bar{k}_j^*$. Finally, define $\bar{S}_j^* := \min_{S_j \in \mathbf{S}_j} \bar{C}_{\bar{k}^* \bar{i}'^* j}(S_j)$, which can be again efficiently solved with the BS algorithm in Appendix A. The following proposition proves that $\left( \bar{i}_j'^*, \bar{k}_j^*, \bar{S}_j^* \right)$ is the optimal solution to sub-problem (26)-(29) with respect to terminal $j$.

**Proposition 1.** $\left( \bar{i}_j'^*, \bar{k}_j^*, \bar{S}_j^* \right) = \min_{i' \in \bar{\mathbf{I}}, k \in \bar{\mathbf{K}}, S_j \in \mathbf{S}_j} \bar{C}_{ki'j}(S_j), \forall j \in \mathbf{J}.$

**Proof.** First, it can be seen from the structure of sub-problem (26)-(29) that as $i'$ varies while the other variables are fixed, $\bar{C}_{ki'j}(S_j)$ increases with $e_{i'j}$. Therefore the optimal solution to $i'$ is $\bar{i}_j'^*$.

Let $\hat{S}_j$ denote the optimal value of $S_j$, then the optimal solution to $k$ is $\hat{k}_j := (\hat{i}_1^j, \hat{i}_2^j, \cdots \hat{i}_L^j) := \min_{k \in \bar{\mathbf{K}}} \bar{C}_{k \bar{i}_j'^* j}(\hat{S}_j)$. We will prove $\hat{k}_j = \bar{k}_j^*$ by contradiction. If there exists $l \in \mathbf{L}$ such that $r_{\hat{i}_l^j j} > r_{\hat{i}_{l+1}^j j}$. We construct a new feasible solution $\tilde{k}_j := (\hat{i}_1^j, \cdots, \hat{i}_{l+1}^j, \hat{i}_l^j, \cdots \hat{i}_L^j)$ by swapping the levels of $\hat{i}_l^j$ and $\hat{i}_{l+1}^j$ in $\hat{k}_j$, and then we compare the difference between the two costs with respect to $\hat{k}_j$ and $\tilde{k}_j$, respectively,

$$\bar{C}_{\hat{k}_j \bar{i}_j'^* j}(\hat{S}_j) - \bar{C}_{\tilde{k}_j \bar{i}_j'^* j}(\hat{S}_j) = (1-q)^2 q^{l-1} d_j \left[ \left( r_{\hat{i}_l^j j} - r_{\hat{i}_{l+1}^j j} \right) + \right.$$

$$e_{\bar{i}_j^{\prime *} j}\left(P_{\hat{i}_l^j j}(\hat{S}_j) - P_{\hat{i}_l^j j}(\hat{S}_j)\right) - \left(r_{\hat{i}_l^j j} P_{\hat{i}_l^j j}(\hat{S}_j) - r_{\hat{i}_{l+1}^j j} P_{\hat{i}_l^j j}(\hat{S}_j)\right)\bigg] >$$

$$(1-q)^2 q^{l-1} d_j \left[\left(r_{\hat{i}_l^j j} - r_{\hat{i}_{l+1}^j j}\right) + \left(e_{\bar{i}_j^{\prime *} j} - r_{\hat{i}_l^j j}\right)\left(P_{\hat{i}_l^j j}(\hat{S}_j) - P_{\hat{i}_l^j j}(\hat{S}_j)\right)\right]$$

Note that $r_{\hat{i}_l^j j} - r_{\hat{i}_{l+1}^j j} > 0$, and $P_{\hat{i}_l^j j}(\hat{S}_j) - P_{\hat{i}_l^j j}(\hat{S}_j) > 0$ due to the assumption that $r_{ij} \geq r_{i'j} \Leftrightarrow t_{ij} \geq t_{i'j}, \forall i \neq i' \in \mathbf{I}, j \in \mathbf{J}$. Then we obtain $\bar{C}_{\hat{k}_j \bar{i}_j^{\prime *} j}(\hat{S}_j) - \bar{C}_{\tilde{k}_j \bar{i}_j^{\prime *} j}(\hat{S}_j) \geq 0$, which is contradictive to the premise that $\hat{k}_j$ is the optimal solution. Therefore, we prove $\hat{i}_1^j \leq \hat{i}_2^j \leq \cdots \leq \hat{i}_L^j$. If there exists a $i \in \bar{\mathbf{I}} \setminus \hat{k}_j$ and some $l \in \mathbf{L}$ having $r_{\hat{i}_l^j j} > r_{ij}$, replacing $\hat{i}_l^j$ with $i$ in $\hat{k}_j$ will further reduce cost $\bar{C}_{\hat{k}_j \bar{i}_j^{\prime *} j}(\hat{S}_j)$ with a similar argument, which is a contradiction, too. This proves that $\hat{k}_j = \bar{k}_j^*$.

Finally, $\hat{S}_j = \bar{S}_j^*$ obviously holds since $\hat{k}_j = \bar{k}_j^*$. This completes the proof. □

By solving problem (26)-(29) for all $j \in \mathbf{J}$, a feasible solution to the primal problem can be obtained as follows:

$$Y_{ijl} = \begin{cases} 1 & \text{if } i = \bar{i}_l^{j*}; \\ 0 & \text{otherwise,} \end{cases} \quad Z_{ij} = \begin{cases} 1 & \text{if } i = \bar{i}_j^{\prime *}; \\ 0 & \text{otherwise,} \end{cases} \quad S_j = \bar{S}_j^*, \quad \forall j \in \mathbf{J}, i \in \mathbf{I}, l \in \mathbf{L}. \quad (35)$$

This algorithm is fast (only taking a solution time of $O\left(|\mathbf{J}| \max\left\{|\mathbf{I}| \ln(|\mathbf{I}|), L \ln(\bar{S}_j)\right\}\right)$). Plugging these feasible solution values into primal objective function (7), we obtain an upper bound to the optimal objective value as well

### 4.3 Updating Lagrangian multipliers

If the upper bound objective obtained from (19) is equal to the lower bound derived in (11), then this bound is the optimal objective value and the corresponding feasible solution is an optimal solution. On the other hand, while the bounds have a considerable gap, multipliers $\lambda$ and $\mu$ will be updated iteratively based to reduce this gap. A subgradient algorithm is used to complete this iterative procedure as described in Appendix B.

### 5. Numerical example

In this section, a series of numerical examples are presented to test the proposed model and provide useful managerial insights based on the datasets provided in Daskin (1995), i.e., a 49-site network involving 48 continental state capital cities and Washington D.C. The numerical algorithms are coded with MATLAB and implemented on a PC with 3.40 GHz CPU and 8 GB RAM. The LR parameters are set as $\tau = 1, \bar{\tau} = 10^{-3}, K = 5, \theta = 1.005$, and $\bar{K} = 60$. We assume that each site has both a candidate supplier and a terminal facility, and the parameters are generated as follows. We set $h_j = h$, $r_{ij} = c_r \delta_{ij}$, and $t_{ij} = c_l \delta_{ij}, \forall i \in \mathbf{I}, j \in \mathbf{J}$, where $h, c_r$ and $c_l$ are constant coefficients and $\delta_{ij}$ is the great-circle distance between sites $i$ and $j$. Each $e_{ij}$ is set to be an independent realization of a uniformly distributed random variable in interval $[1, 1+c_e] \cdot \max a_{i' \in \mathbf{I}} r_{i'j}$, where $c_e \geq 0$ is a constant scalar. In addition, we assume that each $f_i$ is the product of the corresponding city population and a scalar $c_f$, and each $d_j$ is the product of the corresponding state population and a scalar $c_d$. We set $L = 3$ for all the cases.

Firstly we test model (7)-(15). Table 2 summarizes the results of 4 instances on the 49-site network by varying failure probability $q$, where we set $h = 100, c_r = 0.01, c_e = 1, c_f = 0.02, c_d = 10^{-5}$, and $c_l = 10^{-4}$. The optimal gap between the final feasible objective value and the best relaxed objective is denoted by $G$, the solution time is denoted by $T$. The optimal system total inventory and the optimal number of selected suppliers are denoted by $S$ and $N$, respectively. Moreover, define

$$P^E := \frac{\sum_{i \in \mathbf{I}} \sum_{j \in \mathbf{J}} \sum_{l \in \mathbf{L}} d_j (1-q) q^{l-1} P_{ij}(S_j) Y_{ijl}^*}{\sum_{j \in \mathbf{J}} d_j}$$

as the percentage of demand served by the expedited shipments, where $Y_{ijl}^*$ is the best solution to $Y_{ijl}$. Inventory cost, regular shipment cost, marginal expedited shipment, emergency cost, supplier set-up cost, and total system cost are denoted by $C^H$, $C^R$, $C^M$, $C^E$, $C^F$, and $C$, respectively, as defined in equations (1), (2), (4)-(7).

**Table 2** Numerical results for the 49-site network ($c_f$=0.02, $c_d$=$10^{-5}$, $c_l$=$10^{-4}$).

| # | q | T | G(%) | N | $C^F$ | S | $C^H$ | $C^R$ | $C^M$ | $P^E$ (%) | C |
|---|---|---|---|---|---|---|---|---|---|---|---|
| 1 | 0.1 | 3154 | 0.18 | 12 | 21198 | 165 | 16500 | 7323.2 | 4862.3 | 3.73 | 49883.5 |
| 2 | 0.3 | 3367 | 0.22 | 16 | 34858 | 204 | 20400 | 7932.3 | 13357 | 14.9 | 76547.3 |
| 3 | 0.5 | 3571 | 0.41 | 21 | 42724 | 193 | 19300 | 8419.9 | 19478 | 20.6 | 89921.9 |
| 4 | 0.7 | 3653 | 0.47 | 24 | 53379 | 195 | 19500 | 3192.1 | 27354 | 26.1 | 94833.2 |

We note in Table 2 that all the instances are solved with G<1% in one hour. This indicates that our approach can solve problem instances of a realistic size to a near-optimum solution with a reasonable solution efficiency. When $q$ increases, $C^F$, $C^M$, $P^E$ and $C$ significantly increase, while $C^H$ and $C^R$ seem to increase first and then drop. This indicates that as $q$ rises, all the cost components will increase at first, leading to a sharp growth of the total cost. Nevertheless, when $q$ keeps increasing, regular service is unreliable and keeping a higher inventory is no longer an appealing solution. Instead, a higher percentage of expedited shipments are needed to keep the service quality. Meanwhile more suppliers are installed to shorten the shipping distance and offset the increasing expedited shipment cost.

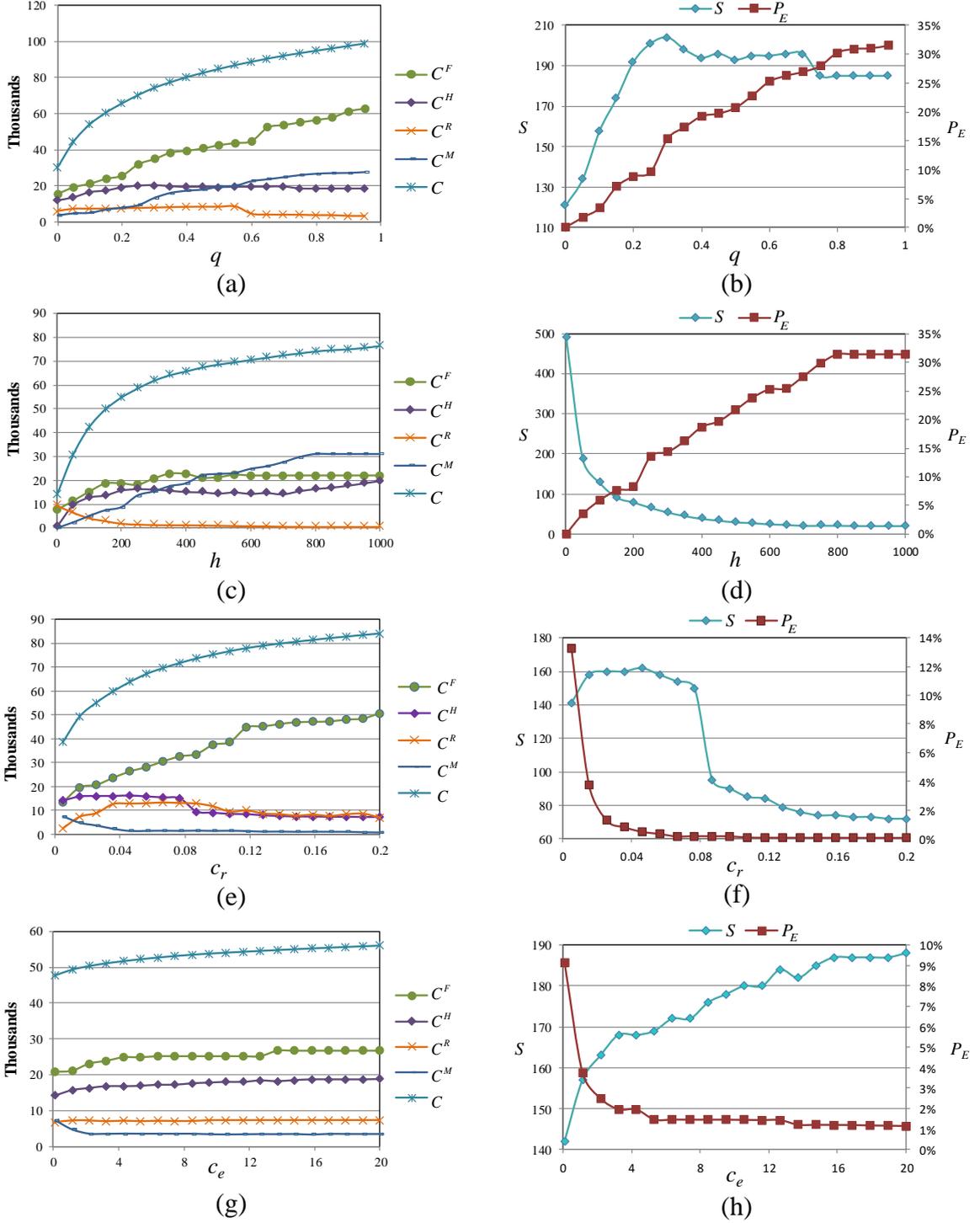

**Figure 2** Sensitivity analysis on parameters $q$ ((a) and (b)), $h$ ((c) and (d)), $c_r$ ((e) and (f)) and $c_e$ ((g) and (h)).

Figure 2 shows four sets of more detailed sensitivity results, where we can see all the cost components, the inventory position $S$ and the expedited service percentage change over key parameters $q$, $h$, $c_r$ and $c_e$. The default parameters are set as $q = 0.1$, $h = 100$, $c_r = 0.01$, $c_e = 1$, $c_f = 0.02$, $c_d = 10^{-5}$, and $c_l = 10^{-4}$, and only one parameter value varies in each experiment. In Figure 2 (a), as $q$ grows from 0 to 1, $C^F$ and $C^M$ generally increase, while $C^H$ increases slightly first and then drops, and $C^R$ is originally stationary and then drops. Also, the total cost $C$ has a sharp increase from around 30000 to 80000, then followed by a constant and slower increasing rate as $q$ becomes lager. Figure 2 (b) shows that $P^E$ rapidly increases with the growth of $q$, while $S$ rises at first and then drops slowly. It's probably because when $q$ increases, the regular service from upstream suppliers become increasingly unreliable, and thus the probability of accessing backup suppliers and expedited services grows. Then more suppliers are selected and higher inventory positions are needed to offset the growth of the shipment costs. Furthermore, as $q$ keeps rising, selecting more suppliers gradually becomes the only cure and higher inventory positions are not as helpful. Meanwhile, expedited shipments gradually take over regular shipments and become the dominating shipment mode. In Figure 2 (c), when $h$ grows from 1 to 1000, $C^F, C^H$ and $C^M$ generally increase, $C^R$ continually drop to almost zero, and $C$ increases strictly first followed by a slower growth. Figure 2 (d) shows that the increase of $h$ rapidly brings down $S$ to a slowing down trend in the tail, while $P^E$ generally increase. This implies that installing more suppliers does not help much when $h$ is large, while using more expedited shipments seems more effective in offsetting the inventory cost growth. We can see in Figure 2 (e) and (f) that both $C$ and $C^F$ increase with the growth of $c_r$ from 0.005 to 0.2, while $P^E$ and $C^M$ keep decreasing to almost zero. $C^H$ and $C^R$ grow slowly at first and then drop, which seems to be consistent with the variation of $S$ in Figure 2(f). This is probably because as $c_r$ grows, the regular shipment cost increases, and thus a higher inventory is needed to offset the growth of expedited shipment cost. The higher inventory leads to a continuous drop of the expedited shipments and a slight increase of the regular shipment cost initially. Nevertheless, as $c_r$ continues to grow (the shipment cost correspondingly increases), building more facilities becomes a better solution to offset the shipment cost growth, which finally brings down the total inventory. In Figure 2(g) and (h), as $c_e$ increases, $S$ increases significantly and $P^E$ drops sharply, but the total cost and all its components do not change too much. This indicates that expedited shipments actually cause little increase in overall cost under the optimal inventory management and transportation configuration strategy, and thus it is an appealing strategy to combine both regular and expedited shipments to reduce the system cost and increase the system reliability.

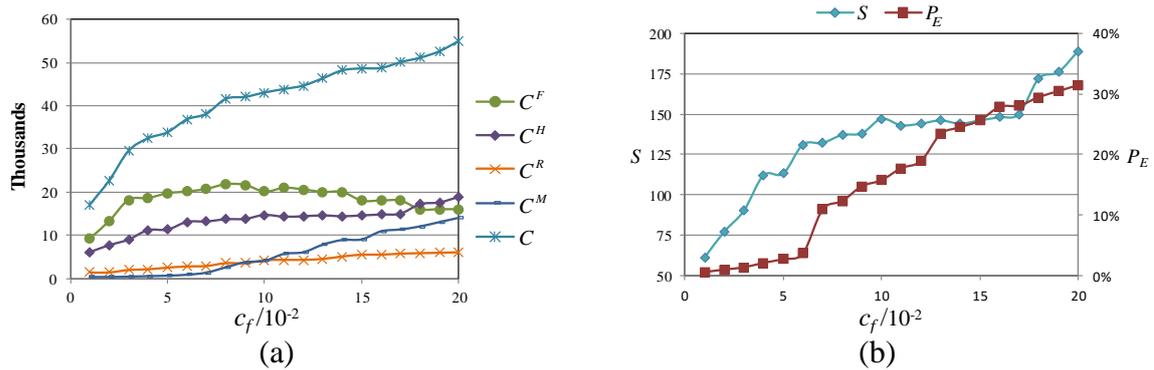

(a)     (b)

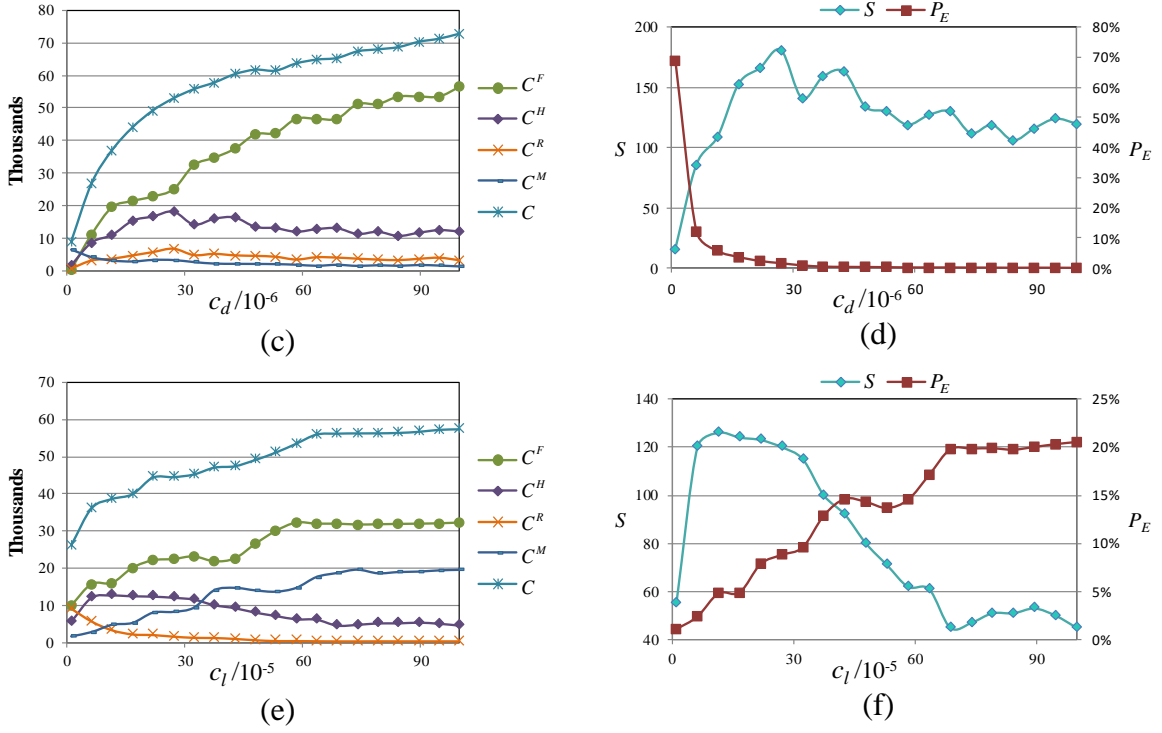

**Figure 3** Sensitivity analysis on parameters $c_f$ ((a) and (b)), $c_d$ ((c) and (d)), and $c_l$ ((e) and (f)).

Besides, we also tested how the results vary with the magnitudes of supplier installation cost (in terms of $c_f$), customer demand rates (in terms of $c_d$) and lead times (in terms of $c_l$), as shown in Figure 3. It can be seen in Figure 3 (a) and (b) that, when $c_f$ initially grows, all cost components and the total cost generally increase. As $c_f$ keeps increasing, $C^F$ increases first and then turns down and $C^R$ flattens out. This is probably because that, the growth of supplier installation cost likely decreases the number of suppliers, which consequentially raises the shipment distance, cost and leading time. Nevertheless, when $c_f$ continues to increase, the number of suppliers is so small that increasing the inventory position alone is not enough to keep the service quality and thus using more expedited services seems necessary. Figure 3 (c) shows that the increase of $c_d$ initially raises all cost components except $C^R$. Then $C^F$ keeps increasing but $C^H$ and $C^R$ decrease slightly with a slowing down trend in the tails. We also see in Figure 3 (d) that $S$ increases quickly initially and then flattens out, while $P^E$ significantly decreases to almost zero. It is probably because that expedited services are more suitable for the cases with low demands when the suppliers are scattered and high inventory positions are unnecessary. Nevertheless, as demands increase, regular shipments will become the main shipment mode instead. Figure 3 (e) and Figure 3 (f) show that as $c_l$ grows, $C^F$, $C^M$ and $P_E$ increase while $C^R$ drop, and $C^H$ and $S$ increases at the beginning then drops. This shows that growth of regular shipment delay will cause the decreasing of inventory positions and consequently increasing the expedited service seems to become a better solution to improve the service quality.

We also test how the variations of $q$ affect the optimal suppliers' layouts and terminal facilities' assignments. Again we set $h = 100, c_r = 0.01, c_e = 1, c_f = 0.02, c_d = 10^{-5},$ and $c_l = 10^{-4}$ and each sub-figure in Figure 4 shows the optimal layout for a different $q$ value among 0, 0.3 and 0.6. In each sub-figure, the squares denote the selected supplier locations and the circles represent the terminal facilities with their area sizes proportional to the base-stock positions. The arrows show how the selected suppliers are assigned to each facility with each arrow's width proportional to the percentage of the corresponding expedited shipments and different colors denoting different levels, i.e., yellow for the first level, green for the second level and pink for the third level.

In Figure 4 (a), when the facility disruption risks are ignored ($q = 0$), the problem would be similar to the integrated model proposed by Li (2013), in which all suppliers are assumed to be reliable and is considered as the benchmark case in our problem. By comparing Figure 4 (a) and (b), we note that as failure probability $q$ increases from 0 to 0.3, 5 more supplier installations are selected and more frequent expedited shipments are needed, in particular for the facilities that are far away from their suppliers. This implies that when primary supplier becomes unreliable and backup facilities are needed, a proper solution is selecting more suppliers to reduce the overall shipment cost. Generally, the expedition percentage increases with the assignment level, and facilities served by farther suppliers hold higher inventory positions.

As $q$ further rises to 0.6, we can see that many more suppliers are installed and the inventory positions of facilities increase generally, but the expedition percentage generally drops. Intuitively, this is because that, as the suppliers become more unreliable, better solutions to offset the growth of shipment cost are increasing the intensity of suppliers distribution and the position of facilities inventory, which will consequently bring down the overall stock-out times and expedition percentage. Also, we can see that 7 more suppliers are installed and most of them (5 suppliers) are located in the northeastern areas with higher population and more facilities. Therefore, under the optimal planning, facilities with more demand may be met first to reduce the shipment costs as much as possible.

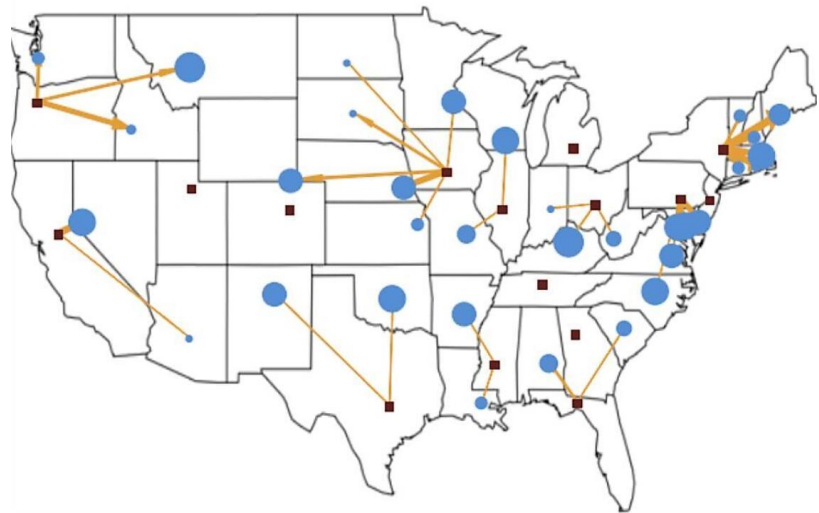

(a) q=0%

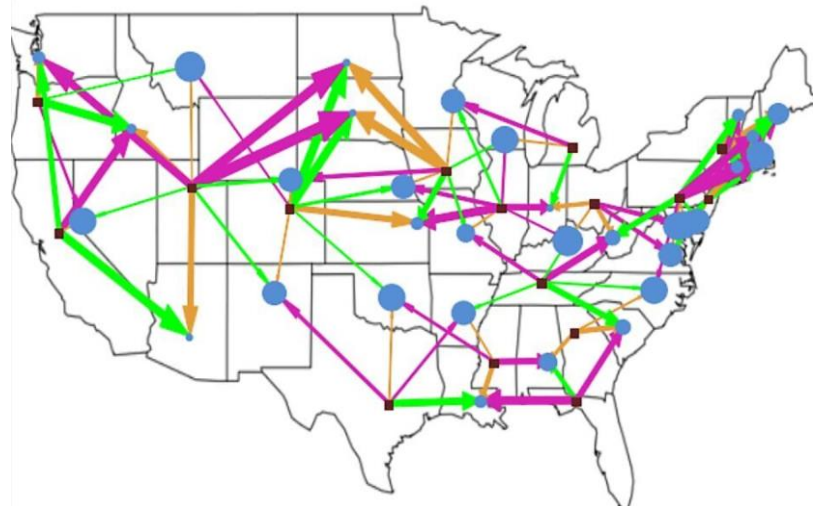
(b) q=30%

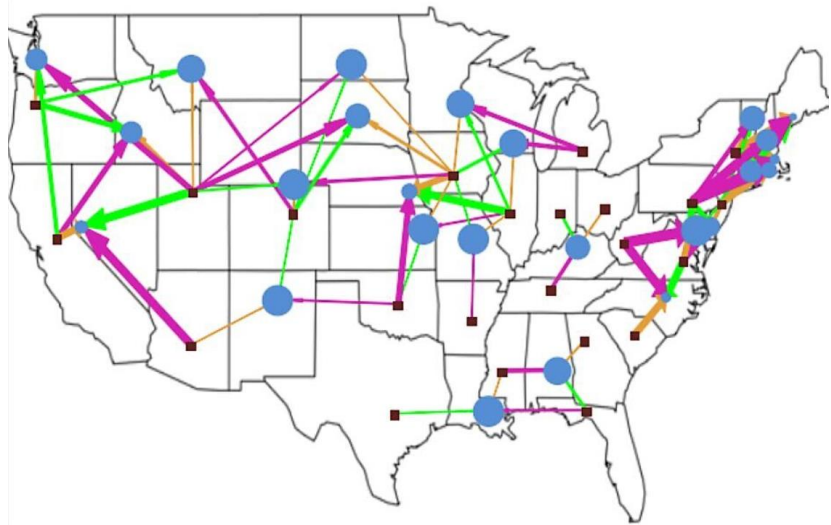
(c) q=60%

**Figure 4** Optimal network layouts for different q.

## 6. Conclusion

This paper presented a reliability model that takes the possibility of supplier failures into the design of integrated logistics system involving logistics network planning and long-term operations management. We formulated a nonlinear mixed integer model to investigate this problem. The model aims to optimize the number of suppliers and their distribution, the assignment of regular shipment services at multiple levels and the expedited services between the suppliers and the terminal facilities, and the base-stock inventory positions. Since the proposed model is difficult to solve (nonlinear integer programming problem), a customized solution approach is created based on Lagrangian relaxation. This solution approach is able to solve this model efficiently and accurately, as evidenced in a set of numerical experiments. Moreover, according to the results of these experiments, a number of managerial insights are found on how the parameter values affect the optimal solution. For example, we found that the optimal network layout and the related cost components can vary significantly with different failure probabilities. Also, when the upstream suppliers become unreliable, expedited shipments will be used more frequently in relation to regular shipments. The increased expedited shipments can ensure the reliability of the integrated logistics system without an excessive increase of the operating costs. The optimal network layout also shows that when the supplier failure probability increases, areas with more customer demands tend to have a high priority to receive services, despite the higher transportation costs.

The proposed model can be further improve in several directions. First, this study assumes that the expedited shipment is "non-fallible", which may be not realistic for some applications where the suppliers may suffer serious disasters causing both services failed. Second, it might be worth considering positive lead times even for expedited shipments for some applications where the expedited delivery time is still noticeable. Finally, this study is set as a two-echelon system where the locations and demand of all the terminal facility are considered to be in the basic conditions and independent of the network design results. However, in some other applications, a more general structure is needed for terminal facilities distribution planning. Extending the current two-echelon network to a more general structure will be an interesting research topic.

**Acknowledgement**

This research is supported by in part by the U.S. National Science Foundation through Grants CMMI #1558889 and CMMI #1541130, the National Center for Intermodal Transportation for Economic Competitiveness (US Tier I University Transportation Center) through Grant # 13091085 and the National Natural Science Foundation of China through Grant #51478151.

**Appendix**

**A Bisecting algorithm to solve (22)**

**Step BS0:** For a given set of $k \in \mathbf{K}, i' \in \mathbf{I}, j \in \mathbf{J}$, initialize two search bounds of as $S_L := 0$ and $S_U := \bar{S}_j$, and the difference slope of $C_{ki'j}(S_j)$ defined in equation (22) with respect to $S_L$ as:

$$G_L := h_j + d_j \sum_{l \in \mathbf{L}} \left(e_{i'j} - r_{i,j}\right)(1-q)q^{l-1}\left(P_{i,j}(0) - P_{i,j}(1)\right)$$

, and that with respect to $S_U$ as:

$$G_U := h_j + d_j \sum_{l \in \mathbf{L}} \left(e_{i'j} - r_{i,j}\right)(1-q)q^{l-1}\left(P_{i,j}(\bar{S}_j) - P_{i,j}(\bar{S}_j - 1)\right).$$

**Step BS1:** If $G_L, G_U \geq 0$, set optimal $S^* := S_L$. Or if $G_L, G_U < 0$, set optimal $S^* := S_U$. Or if $S_U - S_L \leq 1$, set $S^* := S_L$ if $C_{ki'j}(S_L) \leq C_{ki'j}(S_U)$ or $S^* := S_U$ otherwise. If $S^*$ is found, go to Step BS3.

**Step BS2:** Set the middle point $S_M := \lfloor (S_L + S_U)/2 \rfloor$. Calculate the slope at $S_M$ as:

$$G_M := h_j + d_j \sum_{l \in \mathbf{L}} \left(e_{i'j} - r_{i,j}\right)(1-q)q^{l-1}\left(P_{i,j}(0) - P_{i,j}(1)\right) \text{ if } S_M = 0 \text{ or}$$

$$G_M := h_j + d_j \sum_{l \in \mathbf{L}} \left(e_{i'j} - r_{i,j}\right)(1-q)q^{l-1}\left(P_{i,j}(S_M) - P_{i,j}(S_M - 1)\right) \text{ otherwise. If } S_M > 0,$$

set $S_U = S_M$ and $G_U = G_M$; otherwise, set $S_L = S_M$ and $G_L = G_M$. Go to Step BS1.

**Step BS3:** Return $S_j^*$ and $C_{ki'j}(S_j^*)$ as the optimal solution and the optimal objective value to problem (22), respectively.

## B Subgradient algorithm to update Lagrangian multipliers

**Step SG0:** Set initial multipliers $\lambda_{ij}^0 = \mu_{ij}^0 = 0, \forall i \in \mathbf{I}, j \in \mathbf{J}$. Set an auxiliary scalar $0 < \tau \leq 2$ and an iteration index $k := 0$. Set the best known upper bound objective $C := +\infty$.

**Step SG1:** Solve relaxed problem $\Delta(\boldsymbol{\lambda}^k, \boldsymbol{\mu}^k)$ with the solution approach proposed in Section 4.1, and $\{X_i^k\}, \{Y_{ij}^k\}, \{Z_{ij}^k\}, \{S_j^k\}$ denote its optimal solution. If the objective value of $\Delta(\boldsymbol{\lambda}^k, \boldsymbol{\mu}^k)$ does not improve in $K$ consecutive iterations (where $K$ is a predefined number, e.g., 5), we update $\tau = \tau/\theta$, where $h$ is a contraction ratio greater than it.

**Step SG2:** Adapt $\{X_i^k\}, \{Y_{ij}^k\}, \{Z_{ij}^k\}, \{S_j^k\}$ to a set of feasible solution with the algorithm described in Section 4.2. Set $C$ equal to this feasible objective if $C$ is greater than it.

**Step SG3:** Calculate the step size as follows[2]

$$t_k := \frac{\tau(C - \Delta(\boldsymbol{\lambda}^k, \boldsymbol{\mu}^k))}{\sum_{i \in \mathbf{I}, j \in \mathbf{J}} \left( \left( \sum_{l \in \mathbf{L}} Y_{ijl}^k - X_i^k \right)^+ + \left( Z_{ij}^k - X_i^k \right)^+ \right)}.$$

Then update multipliers as follows

$$\lambda_{ij}^{k+1} = \left[ \lambda_{ij}^k + t_k \left( \sum_{l \in \mathbf{L}} Y_{ijl}^k - X_i^k \right) \right]^+, \mu_{ij}^{k+1} = \left[ \mu_{ij}^k + t_k \left( Z_{ij}^k - X_i^k \right) \right]^+, \forall i \in \mathbf{I}, j \in \mathbf{J}.$$

**Step SG4:** Terminate this algorithm if (i) optimality gap $\frac{C - \Delta(\boldsymbol{\lambda}^k, \boldsymbol{\mu}^k)}{C} \leq \varepsilon$ where $\varepsilon$ is a pre-specified error tolerance, (ii) $\tau$ is smaller than a minimum value $\bar{\tau}$, or (iii) $k$ exceeds a maximum iteration number $\bar{K}$; return the best feasible solution as the near-optimum solution. Otherwise $k = k + 1$, and go to Step SG1.

---

[2] In the denominator of this formula, we use the absolute value instead of the squared Euclidean norm, because we found it helps improve the solution efficiency.